\documentclass[11pt]{amsart}
\usepackage{amsmath}

\newtheorem{theorem}{Theorem}[section]

\theoremstyle{definition}

\numberwithin{equation}{section}
\input{tcilatex}

\begin{document}
\title[On Li's coefficients]{Effective method of computing Li's coefficients
and their properties}
\author{Krzysztof Ma\'{s}lanka}
\address{Astronomical Observatory of the Jagiellonian University \\
Orla 171, 30-244 Cracow, Poland}
\date{9 May 2004}
\thanks{e-mail: maslanka@oa.uj.edu.pl}
\keywords{Riemann zeta function, Riemann hypothesis, Li's criterion,
numerical methods in analytic number theory}

\begin{abstract}
In this paper we present an effective method for computing certain real
coefficients $\lambda _{n}$ which appear in a criterion for the Riemann
hypothesis proved by Xian-Jin Li. With the use of this method a sequence of
over three-thousand $\lambda _{n}$'s has been calculated. This sequence
reveals a peculiar and unexpected behavior: it can be split into a strictly
growing \textit{trend} and some tiny \textit{oscillations} superimposed on
this trend.
\end{abstract}

\maketitle

\section{Introduction}

Since its formulation almost a century and a half ago the Riemann hypothesis
(hereafter called RH) is commonly regarded as both the most challenging and
the most difficult task in number theory \cite{Riemann}. It states that all
complex zeroes of the zeta function, defined by the following series if $\Re
s>1$ 
\begin{equation}
\zeta (s)=\sum_{n=1}^{\infty }\frac{1}{n^{s}}  \label{zeta}
\end{equation}%
and by analytic continuation to the whole plane, are located right on the
critical line $\Re s=\frac{1}{2}$. RH, if true, would shed more light on our
knowledge of the distribution of prime numbers. More precisely, the absence
of zeroes of $\zeta (s)$ in the half-plane $\Re s>\theta $ implies that (see 
\cite{Ingham}, theorem 30) 
\begin{equation}
\pi (x)=\mathrm{li}(x)+O(x^{\theta }\log x)  \label{theorem}
\end{equation}%
where $\pi (x)$ is the number of primes not exceeding $x$ and \textrm{li}$%
(x) $ denotes logarithmic integral. Therefore, the value $\theta =\frac{1}{2}
$ (as Riemann conjectured) makes the theorem useful since the error term in (%
\ref{theorem}) is the smallest possible. We do know that on the critical
line lie infinitely many complex zeroes \cite{Hardy} and that among several
billions of initial zeroes there is no counterexample to RH, cf. \cite%
{Odlyzko}, \cite{Wedeniwski}.\newline
\newline
\textit{(Place Figure 1 about here.)}

\section{Li's Criterion}

In 1997 Xian-Jin Li \cite{Li} presented an interesting criterion equivalent
to the Riemann hypothesis:

\begin{theorem}
RH is true if and only if all coefficients%
\begin{equation}
\lambda _{n}:=\frac{1}{\Gamma \left( n\right) }\left. \frac{d^{n}}{ds^{n}}%
\left[ s^{n-1}\ln \xi \left( s\right) \right] \right\vert _{s=1}
\label{lambda1}
\end{equation}%
are non-negative, where 
\begin{equation}
\xi \left( s\right) =2\left( s-1\right) \pi ^{-s/2}\Gamma \left( 1+\frac{s}{2%
}\right) \zeta \left( s\right) .  \label{ksi}
\end{equation}
\end{theorem}

An equivalent definition of $\lambda _{n}$ is (see \cite{Li}, formula 1.4):%
\begin{equation}
\lambda _{n}=\sum_{\rho }\left( 1-\left( 1-\frac{1}{\rho }\right) ^{n}\right)
\label{lambda2}
\end{equation}%
where the sum runs over all (paired) complex zeroes of the Riemann zeta-
function. However, the above definitions of $\lambda _{n}$ are not suitable
for numerical calculations. In this paper I shall present an effective
method for calculating these coefficients. The gathered data investigated
numerically up to $n=3300$ reveals unexpected properties: it contains a
strictly growing trend plus extremely small oscillations superimposed on
this trend.

The following decomposition of $\lambda _{n}$ is implicitly given in a
recent paper by Bombieri and Lagarias (\cite{BombieriLagarias}, Theorem 2):%
\begin{eqnarray}
\lambda _{n} &=&\underset{\text{trend}}{\underbrace{1-\left( \log \left(
4\pi \right) +\gamma \right) \frac{n}{2}+\sum_{j=2}^{n}\left( -1\right) ^{j}%
\binom{n}{j}\left( 1-2^{-j}\right) \zeta \left( j\right) }}
\label{decomposition} \\
&&-\underset{\text{oscillations}}{\underbrace{\sum_{j=1}^{n}\binom{n}{j}\eta
_{j-1}}}  \notag \\
&\equiv &\overset{-}{\lambda }_{n}+\overset{\sim }{\lambda }_{n}  \notag
\end{eqnarray}%
Using the language of signal theory (perhaps not very common but sometimes
appropriate in number theory) one can say that the decomposition (\ref%
{decomposition}) uniquely \textquotedblleft splits" the behavior of the
sequence of $\{\lambda _{n}\}$ into a strictly growing \textit{trend} $%
\overset{-}{\lambda }_{n}$ and certain tiny\textit{\ oscillations} $\overset{%
\sim }{\lambda }_{n}$ superimposed on it. It may be proved that the trend
can be expressed as 
\begin{equation}
\overset{-}{\lambda }_{n}=\frac{1}{\Gamma \left( n\right) }\left. \frac{d^{n}%
}{ds^{n}}\left[ s^{n-1}\ln \left( \pi ^{-s/2}\Gamma \left( 1+\frac{s}{2}%
\right) \right) \right] \right\vert _{s=1}  \label{trend}
\end{equation}%
while the oscillations are%
\begin{equation}
\overset{\sim }{\lambda }_{n}=\frac{1}{\Gamma \left( n\right) }\left. \frac{%
d^{n}}{ds^{n}}\left[ s^{n-1}\ln \left( \left( s-1\right) \zeta \left(
s\right) \right) \right] \right\vert _{s=1}  \label{osc}
\end{equation}

It may also be proved that the trend is indeed strictly growing as $n$ tends
to infinity. It is evident that (\ref{trend}) differs from the main
definition (\ref{lambda1}) simply by replacing $\xi (s)$ by the much simpler
function $\pi ^{-s/2}\Gamma (1+s/2)$. On the other hand, the oscillatory
behavior of $\overset{\sim }{\lambda }_{n}$is not so evident, nevertheless
it may be investigated numerically. I shall return to this decomposition
later.

The problem of calculating numerically both components of $\lambda _{n}$
(i.e. trend and oscillations), using directly (\ref{trend}) and (\ref{osc}),
is rather hopeless, not to say: malicious. One has to take the $n^{th}$
derivatives of functions which depend of variable $s$ and are labelled by
parameter $n$. When $n$ tends to infinity both families of differentiated
functions tend to right angle shaped figures and the derivatives are to be
taken just at the almost singular point $s=1$. What is interesting, the
first derivative in $s=1$ for \textit{all} functions related to the
oscillating part is the same and equal to the Euler constant.

However, it is possible to calculate several tens of initial derivatives
using direct numerical approach, for example Mathametica's ND built-in
function. This function has several parameters which enable to control the
required accuracy. One must be aware, however, that there is no guarantee
that the result will be correct. As with all numerical techniques for
evaluating the infinite via finite samplings, it is sometimes possible to
"fool" ND into giving an incorrect result.\newline
\newline
\textit{(Place Figure 2 about here.)}

It turns out that the following function fits very well to the numerically
tabulated values of (\ref{trend}):%
\begin{equation}
a\left( 1+n\ln n\right) +cn  \label{fit}
\end{equation}%
with%
\begin{eqnarray*}
a &=&\frac{1}{2}\pm 8\cdot 10^{-9} \\
c &=&-1.130330701...
\end{eqnarray*}%
(A choice very similar to (\ref{fit}) was suggested to me by J. Lagarias, 
\cite{LagariasMail}.) Recently I\ learned that A. Voros \cite{Voros}\ using
classic technique of saddle-point method calculated the exact value of $c$%
\begin{equation*}
c=\frac{1}{2}\left( \gamma -1-\ln 2\pi \right)
\end{equation*}%
where $\gamma $ is Euler constant. Simple fitting procedure gives also
coefficients of consecutive terms which appear to be related to Bernoulli
numbers $B_{k}$%
\begin{equation}
-\sum\limits_{k=1}\frac{B_{k}}{2kn^{k-1}}=\frac{1}{4}-\frac{1}{24n}+\frac{1}{%
240n^{3}}-\frac{1}{504n^{5}}+\frac{1}{480n^{7}}-...  \label{asymptotic}
\end{equation}%
This series works well for a dozen or so initial terms although it is
formally divergent.

The starting point of Li's approach to RH is a certain transformation of the
complex plane into itself using the map $s\mapsto z=1-1/s$ (which is a
special case of M\"{o}bius transformation). Under this transformation the
half-plane $\Re s>\frac{1}{2}$ is mapped into the unit disk (with the
critical line $\Re s=\frac{1}{2}$ becoming the unit circle, see Figs. 3 and
4). This was Li's original idea. However, he was inspired by studying A.
Weil's proof of RH for function fields over finite fields where the critical
line is transformed into a unit circle \cite{LiMail}. \newline
\newline
\textit{(Place Figures 3 and 4 about here.)}

\section{The main derivation}

It has been known since medieval times that the harmonic series $%
\sum_{k=1}^\infty 1/k$ diverges. This was proved long ago with the use of
elementary methods by Nicole d'Oresme in the 14th century, and, much later,
independently by Pietro Mengoli (in his book on arithmetic series \textit{%
Novae quadraturae arithmeticae}, 1650) as well as, using yet another method,
by the Bernoulli brothers.

A natural question emerges: how fast does this series diverge? It turns out
that its divergence is \textquotedblleft weak", more precisely: logarithmic.
The quantitative answer to this question implies the definition of the
following famous number called the Euler-Mascheroni constant: 
\begin{equation}
\gamma :=\lim_{x\rightarrow \infty }\left( \sum_{k\leq x}\frac{1}{k}-\log
x\right) =0,5772156649...  \label{gamma}
\end{equation}%
Its natural generalization is the sequence $\gamma _{n}$ defined by 
\begin{equation}
\gamma _{n}:=\frac{(-1)^{n}}{n!}\lim_{x\rightarrow \infty }\left(
\sum_{k\leq x}\frac{1}{k}(\log k)^{n}-\frac{(\log x)^{n+1}}{n+1}\right)
\label{gamman}
\end{equation}%
where $\gamma _{0}=\gamma $. These are the so-called Stieltjes constants%
\footnote{%
It should be noted that the function \textbf{StieltjesGamma[n]} implemented
in Wolfram's \textit{Mathematica},which employs Keiper's algorithm \cite%
{Keiper}, uses a different convention. It is related to our $\gamma _{n}$
via 
\begin{equation*}
\gamma _{n}=\frac{\left( -1\right) ^{n}}{n!}\text{\textbf{StieltjesGamma[n]}}
\end{equation*}%
}. Another \textquotedblleft similar" very useful sequence denoted by $\eta
_{n}$ is defined by 
\begin{equation}
\eta _{n}:=\frac{\left( -1\right) ^{n}}{n!}\lim_{x\rightarrow \infty }\left(
\sum_{k\leq x}\frac{\Lambda \left( k\right) }{k}(\log k)^{n}-\frac{(\log
x)^{n+1}}{n+1}\right) ,  \label{etan}
\end{equation}%
where $\Lambda (k)$ is the so-called von Mangoldt function defined for any
positive integer $k$ as: 
\begin{equation}
\Lambda (k)=\left\{ 
\begin{array}{l}
\log p\text{ if }k\text{ is a prime }p\text{ or any power of a prime }p^{n}
\\ 
0\text{ otherwise}%
\end{array}%
\right.  \label{mangoldt}
\end{equation}

The above sequences are important on their own right since they appear in
the Laurent expansions for $\zeta (s)$ and its logarithmic derivative around 
$s=1$. (There are different conventions when defining these numbers, here I
have adopted those of Bombieri and Lagarias \cite{BombieriLagarias}): 
\begin{equation}
\zeta \left( s+1\right) =\frac{1}{s}+\sum_{n=0}^{\infty }\gamma _{n}s^{n}
\label{series gamma}
\end{equation}%
\begin{equation}
-\frac{\zeta ^{\prime }}{\zeta }\left( s+1\right) =\frac{1}{s}%
+\sum_{n=0}^{\infty }\eta _{n}s^{n}  \label{series eta}
\end{equation}

Integrating the second equation (\ref{series eta}) with respect to $s$,
inserting the result into the first one and equating coefficients in the
appropriate powers of the variable $s$ one can find explicit relations
between the $\gamma _{n}$ and the $\eta _{n}$:

\begin{eqnarray}
\sum_{n=0}^{\infty }\eta _{n}\frac{s^{n+1}}{n+1} &=&-\log \left(
1+\sum_{n=0}^{\infty }\gamma _{n}s^{n+1}\right)  \label{manipulation} \\
\sum_{n=0}^{\infty }\eta _{n}\frac{s^{n+1}}{n+1} &=&\sum_{k=1}^{\infty }%
\frac{\left( -1\right) ^{k}}{k}s^{k}\left( \sum_{n=0}^{\infty }\gamma
_{n}s^{n}\right) ^{k}  \notag
\end{eqnarray}%
Now introduce the coefficients $c_{n}^{(k)}$ defined by%
\begin{equation*}
\sum_{n=0}^{\infty }c_{n}^{\left( k\right) }s^{n}=\left( \sum_{n=0}^{\infty
}\gamma _{n}s^{n}\right) ^{k}
\end{equation*}%
Employing a certain formula from \cite{GradshteynRyzhik} (formula 0.314,
i.e., raising a power series to an arbitrary integral exponent) one can
express the $c$ coefficients by the following recurrence relations: 
\begin{eqnarray}
c_{0}^{\left( k\right) } &=&\gamma ^{k}  \label{c} \\
c_{m}^{(k)} &=&\frac{1}{m\gamma }\sum_{i=0}^{m-1}\left[ km-\left( k+1\right)
i\right] \gamma _{m-i}c_{i}^{(k)}  \notag
\end{eqnarray}

The matrix of coefficients $c$ depends on $\{\gamma _{n}\}$: 
\begin{equation}
\begin{array}{c}
\bigskip \\ 
k=1 \\ 
k=2 \\ 
k=3 \\ 
k=4 \\ 
k=5%
\end{array}%
\hspace{0.4cm}%
\begin{array}{ccccc}
m=0\bigskip & m=1 & m=2 & m=3 & m=4 \\ 
\gamma _{0} & \gamma _{1} & \gamma _{2} & \gamma _{3} & \gamma _{4} \\ 
\gamma _{0}^{2} & 2\gamma _{0}\gamma _{1} & \gamma _{1}^{2}+2\gamma
_{0}\gamma _{2} & 2\gamma _{1}\gamma _{2}+2\gamma _{0}\gamma _{3} & ... \\ 
\gamma _{0}^{3} & 3\gamma _{0}^{2}\gamma _{1} & 3\gamma _{0}\gamma
_{1}^{2}+3\gamma _{0}^{2}\gamma _{2} & ... & ... \\ 
\gamma _{0}^{4} & 4\gamma _{0}^{3}\gamma _{1} & ... & ... & ... \\ 
\gamma _{0}^{5} & ... & ... & ... & ...%
\end{array}
\label{cmk}
\end{equation}%
(In what follows only the upper triangular part of this infinite matrix will
be needed.)\newline
\newline
\textit{(Place Figure 5 about here.)} \newline
\newline

With the help of (\ref{manipulation}) the coefficients $\eta _{n}$ may
further be expressed using the elements of the matrix $c$ as 
\begin{equation}
\eta _{n}=(n+1)\sum_{k=0}^{n}\frac{(-1)^{k+1}}{k+1}c_{n-k}^{(k+1)}
\label{etan1}
\end{equation}

From this we have: 
\begin{eqnarray}
\eta _{0} &=&-\gamma _{0}  \label{etas} \\
\eta _{1} &=&+\gamma _{0}^{2}-2\gamma _{1}  \notag \\
\eta _{2} &=&-\gamma _{0}^{3}+3\gamma _{0}\gamma _{1}-3\gamma _{2}  \notag \\
\eta _{3} &=&+\gamma _{0}^{4}-4\gamma _{0}^{2}\gamma _{1}+2\gamma
_{1}^{2}+4\gamma _{0}\gamma _{2}-4\gamma _{3}  \notag \\
\eta _{4} &=&-\gamma _{0}^{5}+5\gamma _{0}^{3}\gamma _{1}-5\gamma _{0}\gamma
_{1}^{2}-5\gamma _{0}^{2}\gamma _{2}+5\gamma _{1}\gamma _{2}+5\gamma
_{0}\gamma _{3}-5\gamma _{4},  \notag \\
&&...  \notag
\end{eqnarray}%
Finally, the oscillating parts of $\lambda _{n}$ are expressible as
polynomials in the Stieltjes constants: 
\begin{equation}
\overset{\sim }{\lambda }_{n}=-\sum_{j=1}^{n}\binom{n}{j}\eta _{j-1}.
\label{oscillations}
\end{equation}%
Using now (\ref{etas}) and (\ref{oscillations}) we finally obtain:

\begin{eqnarray}
\overset{\sim }{\lambda }_{1} &=&\gamma _{0}  \label{lambdas} \\
\overset{\sim }{\lambda }_{2} &=&2\gamma _{0}-\gamma _{0}^{2}+2\gamma _{1} 
\notag \\
\overset{\sim }{\lambda }_{3} &=&3\gamma _{0}-3\gamma _{0}^{2}+\gamma
_{0}^{3}+6\gamma _{1}-3\gamma _{0}\gamma _{1}+3\gamma _{2}  \notag \\
\overset{\sim }{\lambda }_{4} &=&4\gamma _{0}-6\gamma _{0}^{2}+4\gamma
_{0}^{3}-\gamma _{0}^{4}+12\gamma _{1}-12\gamma _{0}\gamma _{1}+4\gamma
_{0}^{2}\gamma _{1}-2\gamma _{1}^{2}  \notag \\
&&+12\gamma _{2}-4\gamma _{0}\gamma _{2}+4\gamma _{3}  \notag \\
&&......  \notag
\end{eqnarray}%
Here are some numerical values of various numbers used in this paper:

\begin{tabular}{|l|l|l|l|l|}
\hline
$n$ & $\gamma _{n}$ & $\eta _{n}$ & $\overset{\sim }{\lambda }_{n}$ & $%
\lambda _{n}$ \\ \hline
{\tiny 0} & {\tiny +0.577215664902} & {\tiny -0.577215664902} & {\tiny -} & 
{\tiny -} \\ \hline
{\tiny 1} & {\tiny -0.0728158454837} & {\tiny +0.187546232840} & {\tiny %
0.577215664902} & {\tiny 0.0230957089661} \\ \hline
{\tiny 2} & {\tiny -0.00969036319287} & {\tiny -0.0516886320332} & {\tiny %
0.966885096963} & {\tiny 0.0923457352280} \\ \hline
{\tiny 3} & {\tiny +0.00205383442030} & {\tiny +0.0147516588255} & {\tiny %
1.22069692822} & {\tiny 0.207638920554} \\ \hline
{\tiny 4} & {\tiny +0.00232537006547} & {\tiny -0.00452447788850} & {\tiny %
1.37558813187} & {\tiny 0.368790479492} \\ \hline
{\tiny 5} & {\tiny +0.000793323817301} & {\tiny +0.00144679520453} & {\tiny %
1.45826850020} & {\tiny 0.575542714461} \\ \hline
{\tiny 6} & {\tiny -0.000238769345430} & {\tiny -0.000471544078185} & {\tiny %
1.48829832721} & {\tiny 0.827566012282} \\ \hline
{\tiny 7} & {\tiny -0.000527289567058} & {\tiny +0.000155180294164} & {\tiny %
1.48019084024} & {\tiny 1.12446011757} \\ \hline
{\tiny 8} & {\tiny -0.000352123353803} & {\tiny -0.0000513452121181} & 
{\tiny 1.44485574412} & {\tiny 1.46575567715} \\ \hline
{\tiny 9} & {\tiny -0.0000343947744181} & {\tiny +0.0000170413570471} & 
{\tiny 1.39059640679} & {\tiny 1.85091604838} \\ \hline
{\tiny 10} & {\tiny +0.000205332814909} & {\tiny -5.66605092104$\cdot
10^{-6} $} & {\tiny 1.32380368370} & {\tiny 2.27933936319} \\ \hline
{\tiny 11} & {\tiny +0.000270184439544} & {\tiny +1.88584861186$\cdot
10^{-6} $} & {\tiny 1.24944277582} & {\tiny 2.75036083822} \\ \hline
{\tiny 12} & {\tiny +0.000167272912105} & {\tiny -6.28055422786$\cdot
10^{-7} $} & {\tiny 1.17139824694} & {\tiny 3.26325532062} \\ \hline
{\tiny 13} & {\tiny -0.0000274638066038} & {\tiny +2.09240519074$\cdot
10^{-7}$} & {\tiny 1.09272131711} & {\tiny 3.81724005785} \\ \hline
{\tiny 14} & {\tiny -0.000209209262059} & {\tiny -6.97247031237$\cdot
10^{-8} $} & {\tiny 1.01580941259} & {\tiny 4.41147767868} \\ \hline
{\tiny 15} & {\tiny -0.000283468655320} & {\tiny +2.32371573798$\cdot
10^{-8} $} & {\tiny 0.942538421086} & {\tiny 5.04507937203} \\ \hline
{\tiny 100} & {\tiny -4.25340157171$\cdot 10^{17}$} & {\tiny -6.46775072494$%
\cdot 10^{-49}$} & {\tiny 0.628752815248} & {\tiny 118.603775377} \\ \hline
{\tiny 500} & {\tiny -1.16550527223$\cdot 10^{204}$} & {\tiny -9.16750985401$%
\cdot 10^{-240}$} & {\tiny 2.66350209695} & {\tiny 991.900092992} \\ \hline
{\tiny 1000} & {\tiny -1.57095384420$\cdot 10^{486}$} & {\tiny -2.52129710770%
$\cdot 10^{-478}$} & {\tiny 1.75626461597} & {\tiny 2326.05316169} \\ \hline
{\tiny 2000} & {\tiny +2.68042467892 $\cdot 10^{1109}$} & {\tiny %
-1.90708173159$\cdot 10^{-955}$} & {\tiny 10.7685011806} & {\tiny %
5351.75953838} \\ \hline
{\tiny 3000} &  &  & {\tiny -2.09002802367} & {\tiny 8617.21920730} \\ \hline
\end{tabular}

\section{Applications and conclusions}

The recurrence formulae (\ref{c}) together with (\ref{etan1}) and (\ref%
{oscillations}) allow in principle to compute both $\eta _{n}$ and $\overset{%
\sim }{\lambda }_{n}$ with arbitrary accuracy for any value of $n$, but it
is clear that with increasing $n$ the number of terms increases very rapidly%
\footnote{%
Using the On-Line Encyclopedia of Integer Sequences
(http://www.research.att.com/\symbol{126}njas/sequences/) one can see that
number of terms in $\eta _{n}$ is related to the number of partitions of $n$
(partition number, \textbf{PartitionsP[n]} in \textit{Mathematica}
notation), which grows like $\exp ($const$\sqrt{n})$, whereas the number of
terms in $\overset{\sim }{\lambda }_{n}$ is equal to the number of sums $S$
of positive integers satisfying $S\leq n$ (\textbf{Sum[PartitionsP[k],%
\{k,1,n\}] }in \textit{Mathematica} notation).}. It would be desirable to
simplify the polynomials in (\ref{etas}) and (\ref{lambdas}), or at least to
reveal some hidden regularities in them, but I doubt whether this is
possible. The table below demonstrates that it would be even impractical to
write down explicit expressions for, say, $\overset{\sim }{\lambda }_{n}$
for $n$ greater than 15 or 20.

\begin{tabular}{|c|c|c|}
\hline
$n$ & \# of terms in $\eta _{n}$ (eq. \ref{etas}) & \# of terms in $\overset{%
\sim }{\lambda }_{n}$ (eq. \ref{lambdas}) \\ \hline
$0$ & $1$ & $-$ \\ \hline
$1$ & $2$ & $1$ \\ \hline
$2$ & $3$ & $3$ \\ \hline
$3$ & $5$ & $6$ \\ \hline
$4$ & $7$ & $11$ \\ \hline
$5$ & $11$ & $18$ \\ \hline
$10$ & $56$ & $138$ \\ \hline
$20$ & $792$ & $2713$ \\ \hline
$30$ & $6842$ & $28628$ \\ \hline
$40$ & $44583$ & $215307$ \\ \hline
$50$ & $239943$ & $1295970$ \\ \hline
\end{tabular}%
\vspace{0.4cm}

Using the above formulae (\ref{c}), (\ref{etan1}) and (\ref{lambdas}) I have
computed quite a lot of initial values of $\eta _{n}$ and $\lambda _{n}$.
First it was necessary to tabulate Stieltjes constants $\gamma _{n}$ with
sufficient number of significant digits. In order to obtain these I used 
\textit{Mathematica 5} which can handle arbitrary precision numbers and
performs automatically full control of accuracy in numerical calculations.
(For details concerning \textit{Mathematica} interval arithmetic see e.g. 
\cite{Mathematica}). This part of computations took over 60 hours on AMD
1667 MHz processor.

Recently Kreminski \cite{Kreminski} published an effective method of
computing Stieltjes gamma using Newton-Cotes integration algorithm. His
method would be of considerable interest since, due to some "hardcoded"
limitations, the current \textit{Mathematica} version can't give $\gamma
_{n} $ (with sufficient accuracy) beyond $n\approx 2050$.

The main calculations ($\eta _{n}$ and $\lambda _{n}$) were also time
consuming (about 20 hours) and required also considerable amount of computer
memory (3$\times $256 Mb). In particular, having 2000 pre-computed Stieltjes
constants, with 800 significant digits each, I calculated 2000 $\eta _{n}$
and almost 3300 $\lambda _{n}$. Due to finite accuracy and the obvious
phenomenon of error accumulation, the number of significant digits in $\eta
_{n}$ and $\lambda _{n}$ decreases with increasing $n$ (see Fig. 6).

In order to verify the computations as well as to compare various packages I
also tried to repeat the whole procedure using \textit{Maple 8}. However, it
turned out that it is impossible to get Stieltjes constants $\gamma _{n}$
for relatively small $n$ =100 with the required precision 800 significant
digits.\newline
\newline
\textit{(Place Figure 6 about here.)} \newline
\newline

The main conclusion which stems from the above calculations is contained in
the following plots showing the trend of $\lambda $ (Fig. 7a) and the
oscillating part of $\lambda $ (Fig. 7b). Their sum gives the coefficients
which appear in Li's criterion for RH. Note that the scales on both plots
differ by nearly two orders of magnitude. As mentioned before, it is easy to
show that the trend (\ref{trend}) is strictly growing. Therefore, if the
oscillations were bounded or, at least, if their amplitude would grow with $%
n $ slower than the trend, then RH would be true. In other words, we have a
new RH criterion, which is simply a reformulation the original Li's result,
but from the viewpoint of the present paper it has an obvious
interpretation. It states that if for all positive integer $n$%
\begin{equation*}
-\overset{\sim }{\lambda }_{n}\leq \overset{-}{\lambda }_{n}
\end{equation*}%
then RH\ is true. The numerical data gathered so far and presented in
Figures 7a and 7b is in its favor. Of course, one should bear in mind that
in number theory the numerical evidence, no matter how "convincing", may be
just illusory. In fact, Oesterl\'{e} observed that if the first $n$ zeta
zeros are on the critical line, then the Li positivity should hold for about
the first $n^{2}$ Li coefficients (see \cite{Biane}, p. 441). Therefore,
direct numerical search for a possible counterexample to RH using Li's
criterion is rather a hopeless task.

Finally I would like to stress out that so far there are no published
extensive tables of Li's coefficients. Several numerical values of $\lambda
_{n}$ are given in \cite{Biane}. All the numerical data obtained during
preparation of this paper as well as appropriate \textit{Mathematica}
notebook are available from the author. \newline
\newline
\textit{(Place Figures 7a and 7b about here.)} \newline
\newline
\textbf{Acknowledgments.} I would like to express my gratitude to Prof.
Jeffrey C. Lagarias, AT\&T Labs, for confirming the effect of tiny
oscillations, as well as for several remarks concerning further development
of the idea presented in this paper. I would also like to thank Dr. Mark W.
Coffey, Department of Physics, Colorado School of Mines, for directing my
attention to some misprints as well as for sending me a preliminary version
of his paper prior to publication \cite{Coffey}. Finally, the Referee of
this paper made many remarks as well as valuable suggestions concerning
future development of the ideas presented in this paper.\newline
\textbf{Note added in proof.} After completing the calculations I become
aware of the paper by Keiper \cite{Keiper}. His Figure 1 looks similar to
Figure 7b of the present paper. However it is not exactly the same, insofar
as he is subtracting off an approximation to the "main term". But presumably
the error is a constant plus $O(1/n)$ term, in view of my asymptotic series
approximation (\ref{asymptotic}), so it must be pretty close.\newline

\ \newline
{\LARGE Figure captions}

\begin{itemize}
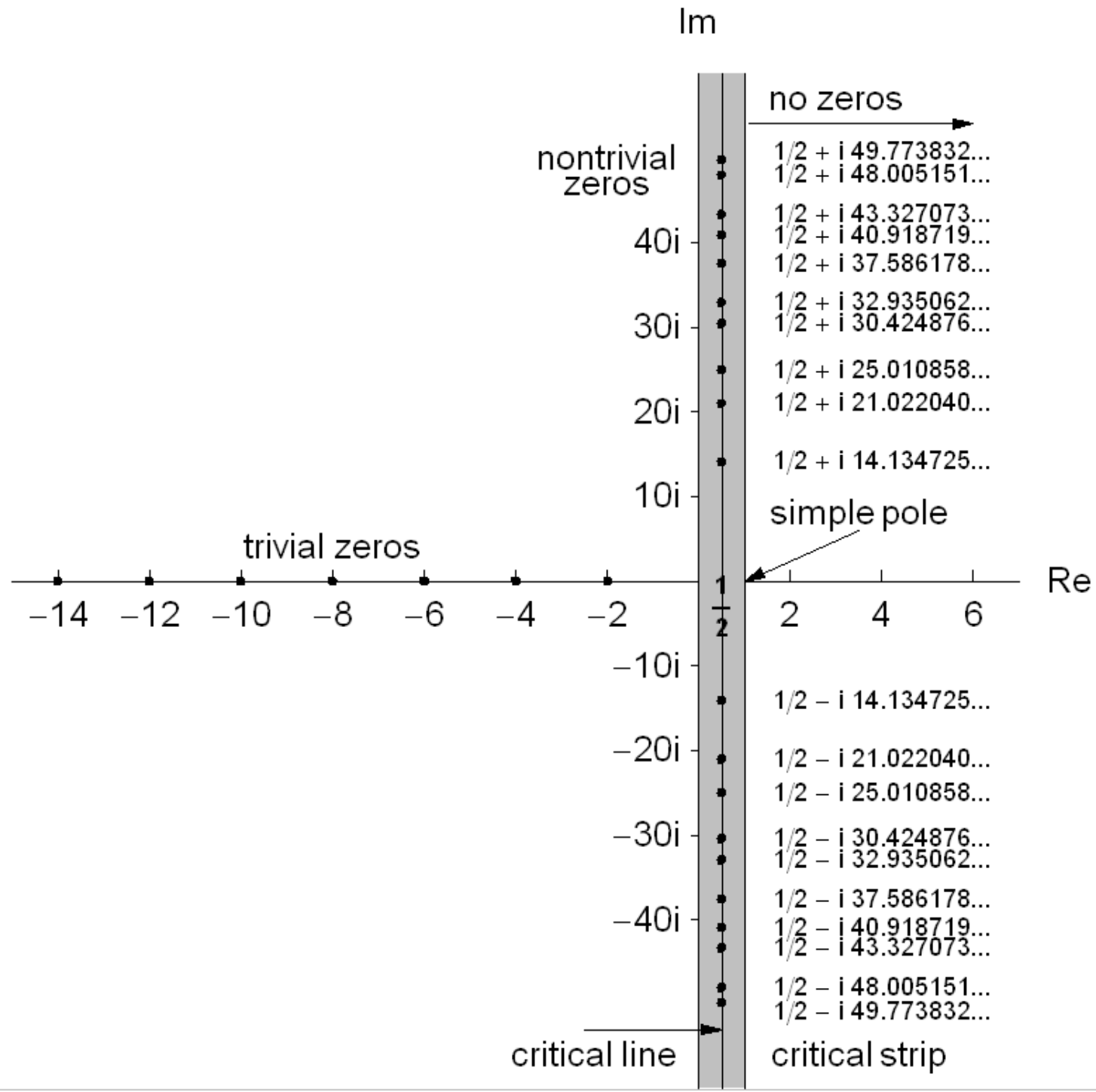
\item Figure 1. Distribution of zeroes of $\zeta $ in the complex plane.

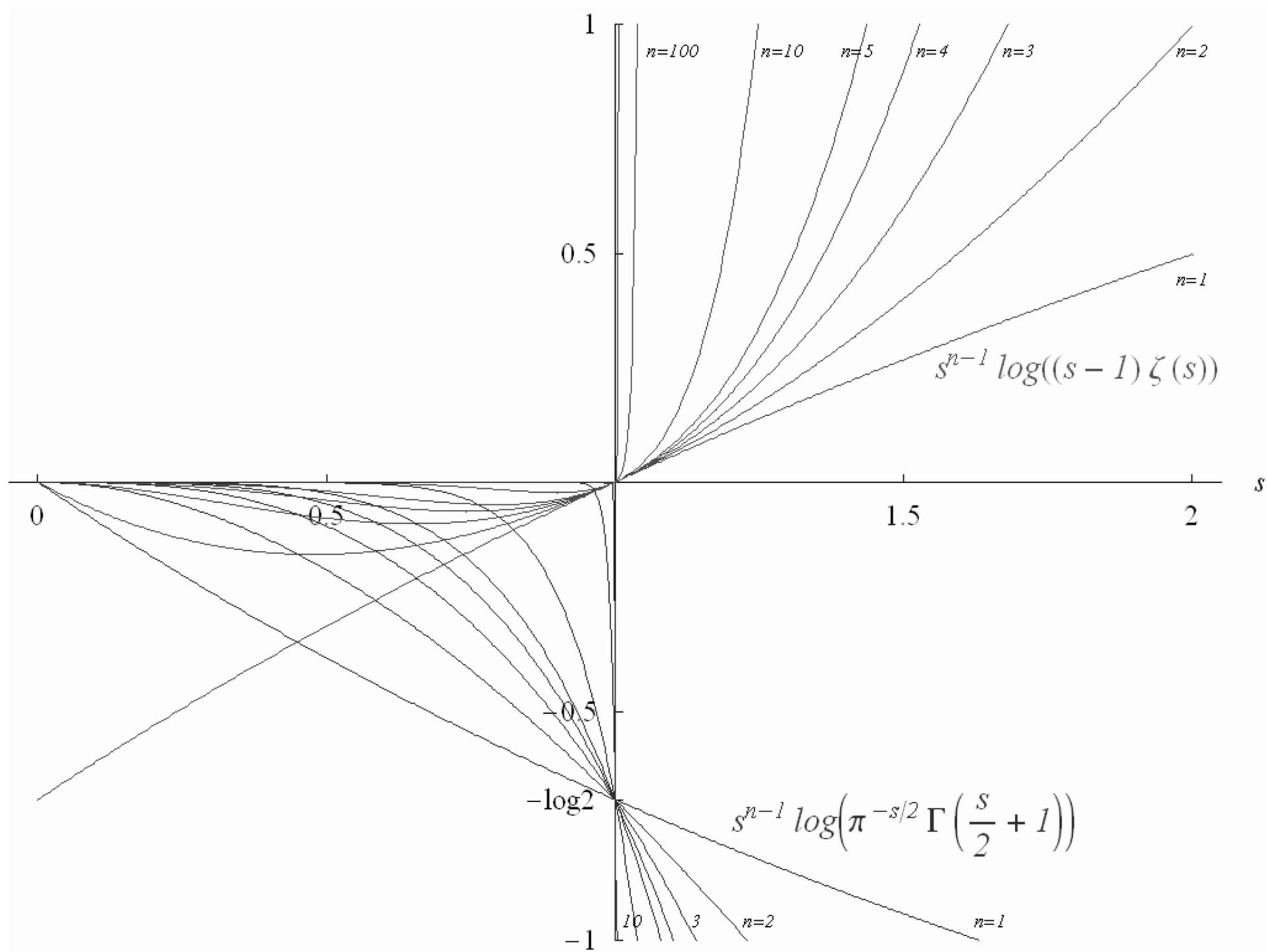
\item Figure 2. Understanding the Li's lambda: inevitable difficulties
encountered when calculating numerically both parts of $\lambda _{n}$. The
sequence of consecutive derivatives is to be taken near $s=0$.

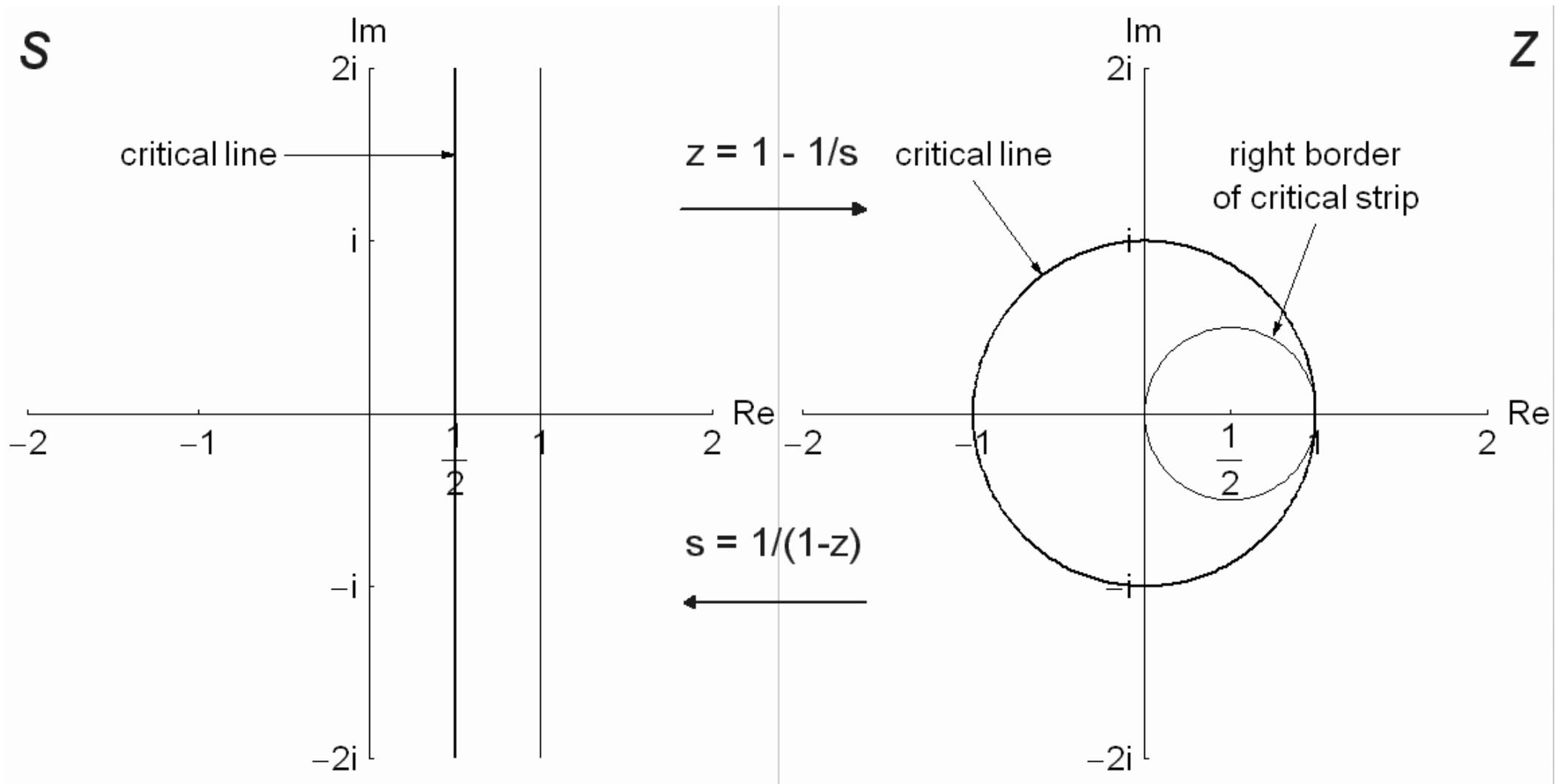
\item Figure 3. M\"{o}bius transformation of the complex plane used by Li.
The left fragment of the picture is just Fig. 1 reduced to its essential
part: the critical strip. The half-plane $\Re s>\frac{1}{2}$ (left picture)
is mapped into the unit disk $|z|<1$ (right picture).

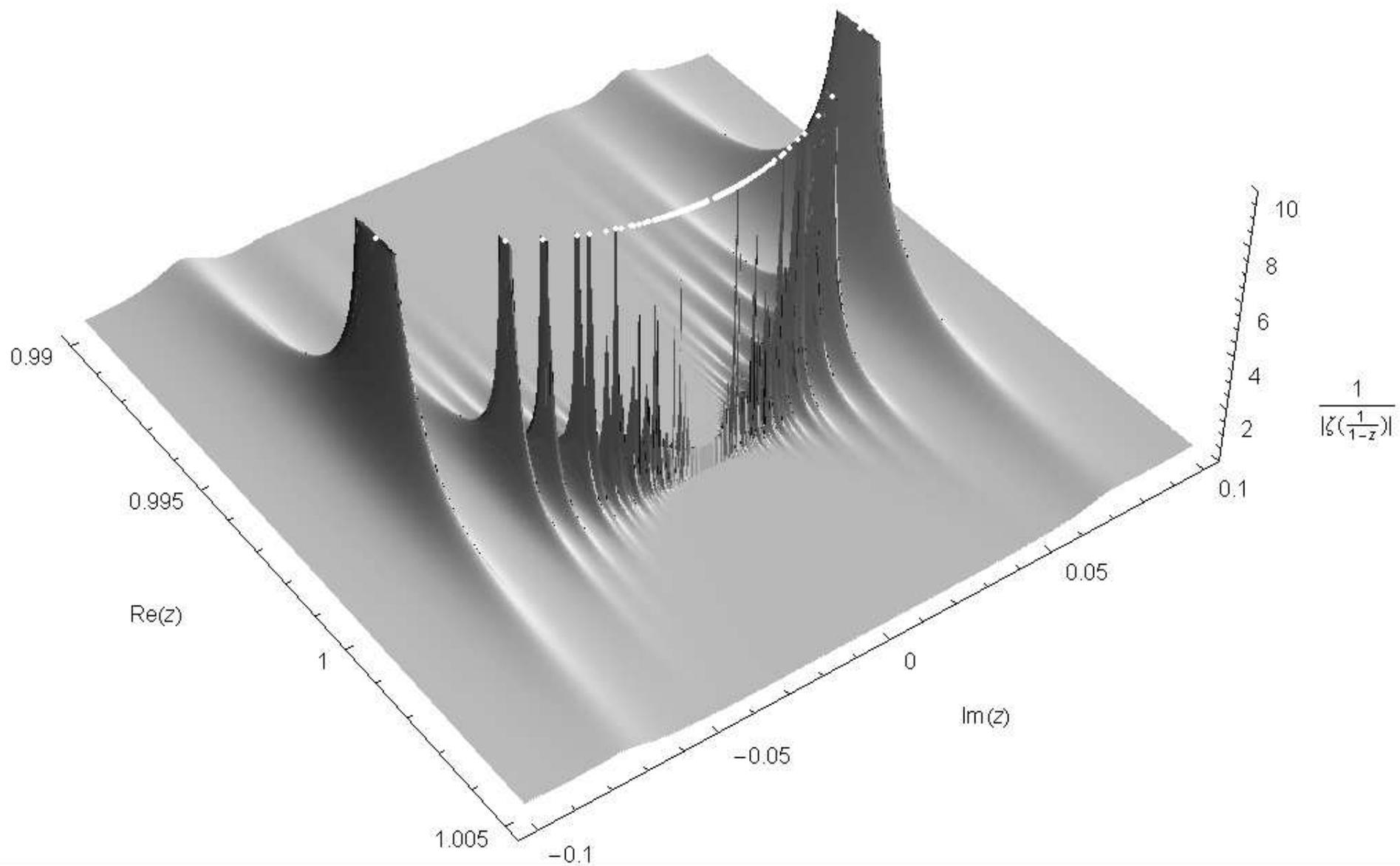
\item Figure 4. Plot of $1/|\zeta (\frac{1}{1-z})|$ on a small part of the
transformed complex plane containing all nontrivial zeroes. On the right
fragment of Fig. 2 this part of complex plane is a small, very narrow
rectangle near $z=1$. Nontrivial zeroes are visible as sharp
\textquotedblleft pins". White dots are added to help visualize that the
zeroes indeed lie on a circle (which looks rather like an ellipse here since
the scales on $\Re z$ and $\Im z$ are different). The apparent lack of peaks
in the center is an artifact. All complex zeroes are very crowded near $z=1$
and the corresponding peaks are increasingly thinner. Obtaining a better
picture would require much higher density of points in which values of zeta
are calculated (hence more computer memory) and much higher resolution of
the picture.

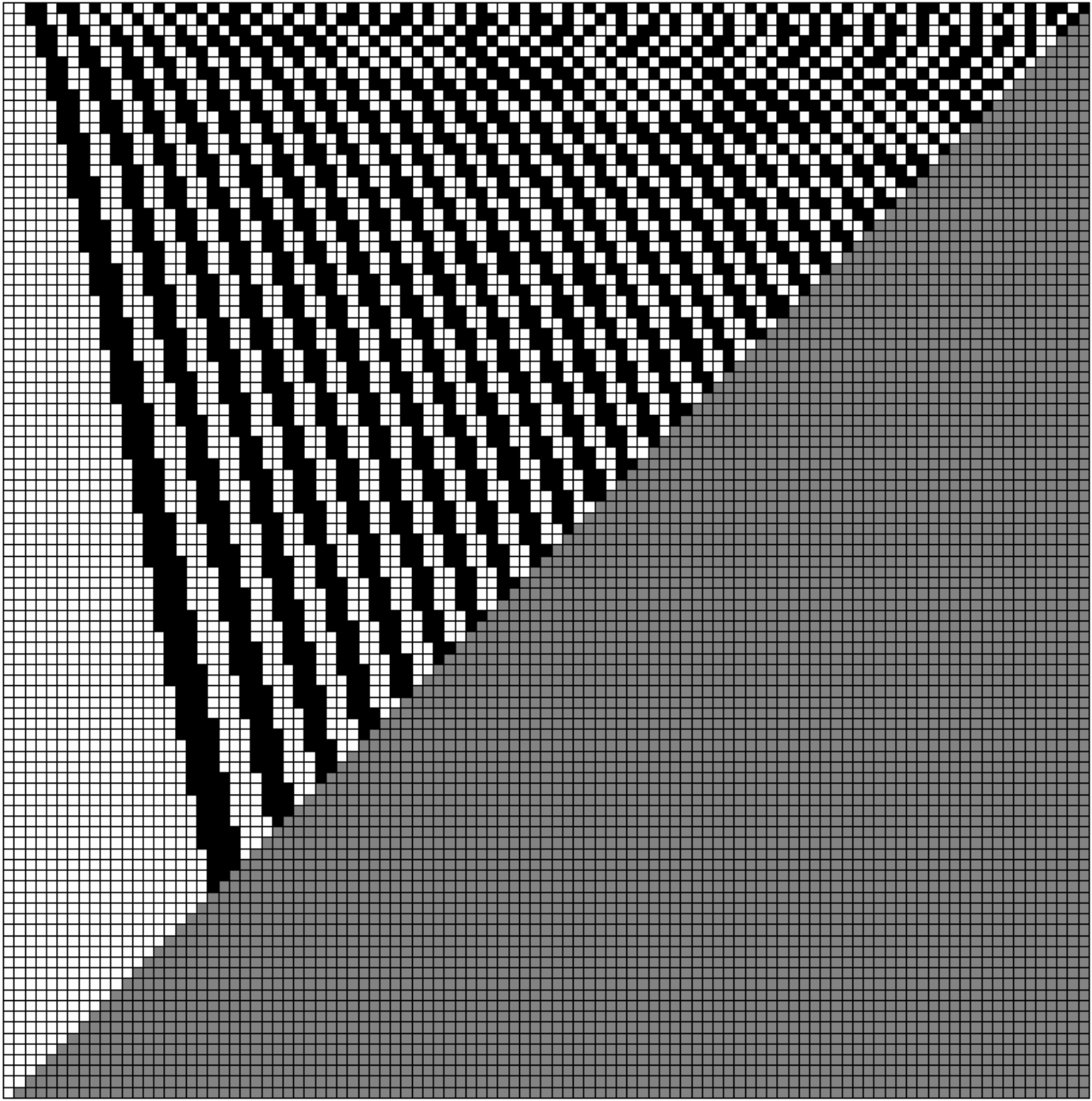
\item Figure 5. Signs of the coefficients of matrix $c$ (\ref{c}) for $k=100$
with rows and columns labelled as in (\ref{cmk}). Little white squares
denote plus sign, black squares denote minus sign; grey squares mark unused
entries of the matrix.

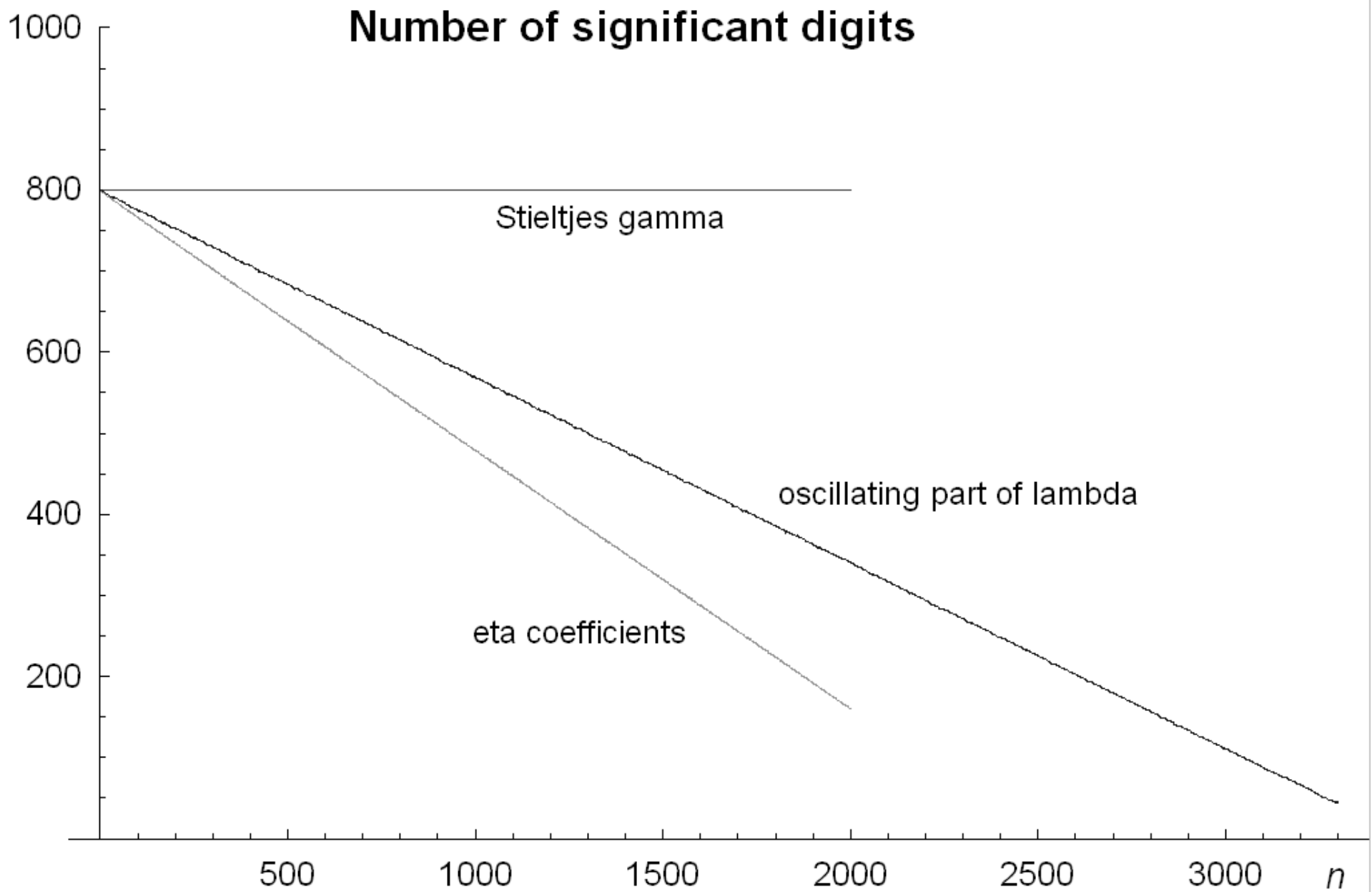
\item Figure 6. Accuracies of various numbers used in this paper: $\gamma
_{n}$, $\eta _{n}$ and $\overset{\sim }{\lambda }_{n}$. Having 2000
precomputed Stieltjes constants $\gamma _{n}$, with 800 significant digits
each, I could (using \textit{Mathematica 5}) obtain 2000 coefficients $\eta
_{n}$ and about 3300 oscillating parts of lambda, $\overset{\sim }{\lambda }%
_{n}$, both with linearly decreasing accuracy. The accuracies of $\eta _{n}$
and $\overset{\sim }{\lambda }_{n}$decrease with $n$ due to complicated
error cumulating but almost perfect linearity of their dependence is a
priori not so obvious because the number of terms in (\ref{etas}) and (\ref%
{lambdas}) grows fast with increasing $n$. In particular, the fact that the
accuracy of $\eta _{n}$ decreases \textit{faster} than accuracy of $\overset{%
\sim }{\lambda }_{n}$ is rather counter intuitive.

\item Figures 7a and 7b. The trend of $\lambda _{n}$ (7a) in comparison with
the oscillating part of $\lambda _{n}$ (7b). Note different vertical scales.
In fact, the sum of the trend and the oscillating part, i.e. full $\lambda
_{n}$, would look exactly like the upper plot since the amplitude of the
oscillations is smaller than the thickness of the graph line.
\end{itemize}

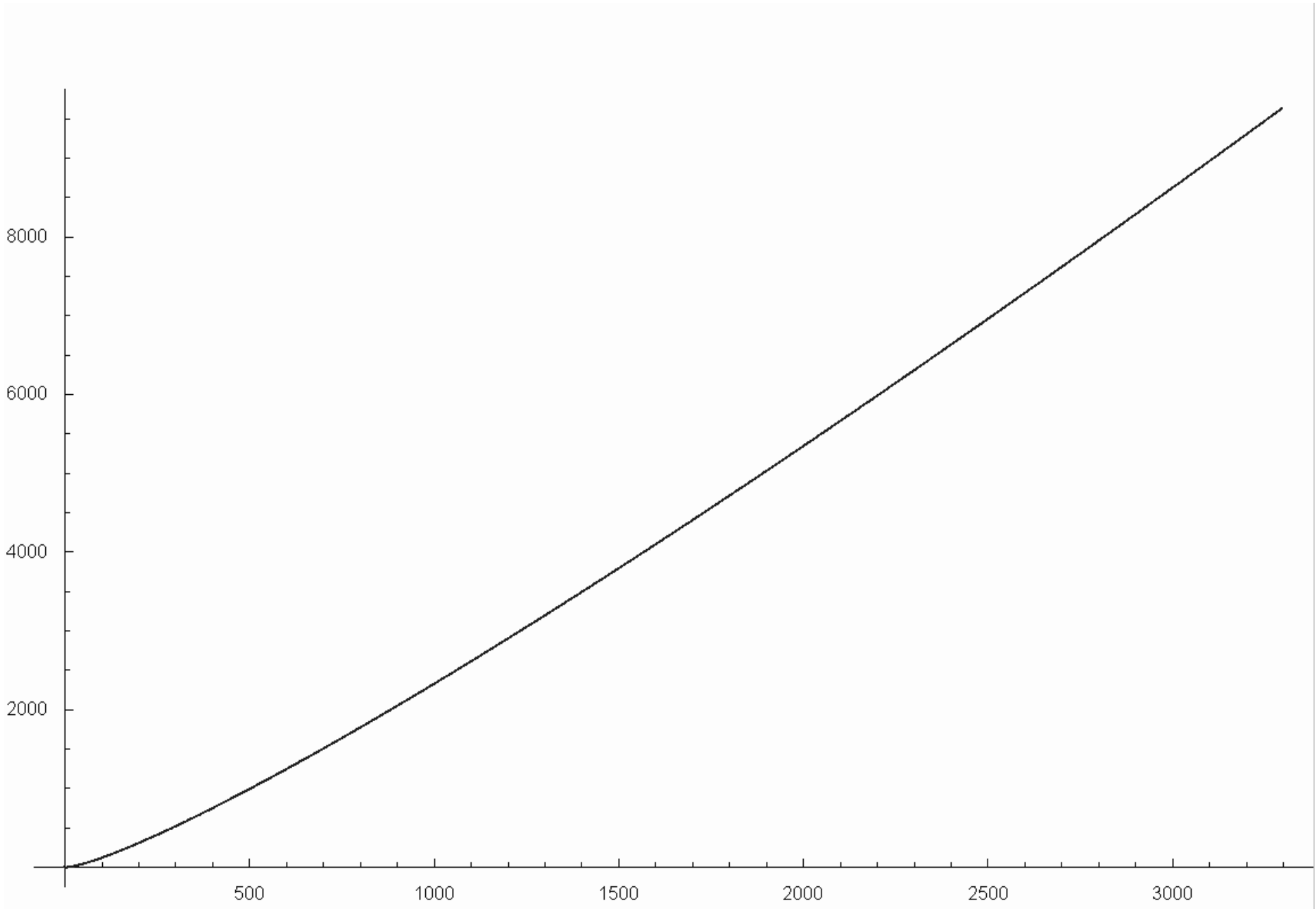
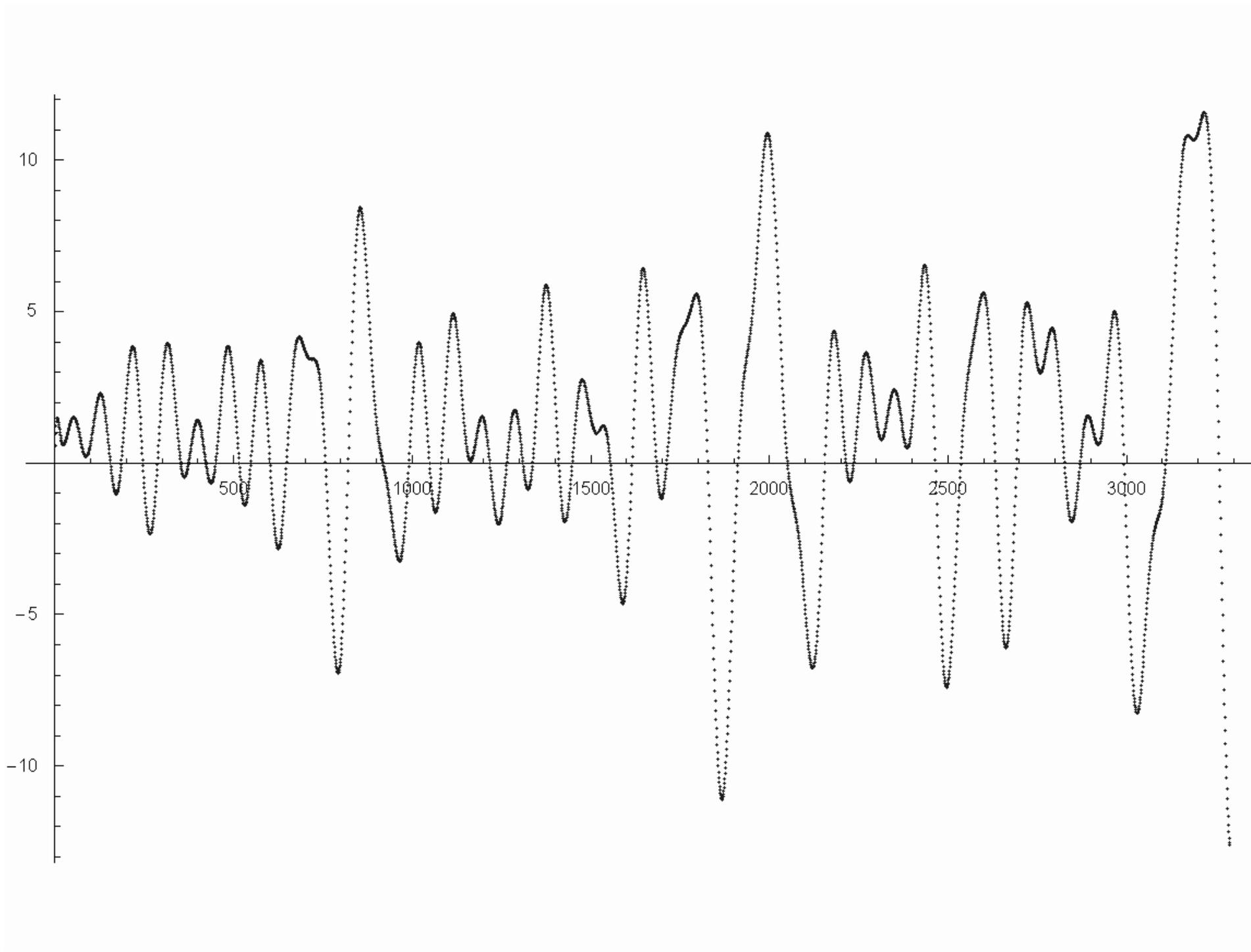
\end{document}